\documentclass[11pt,a4paper,leqno,twoside, headinclude]{amsart}

\usepackage[tracking=true]{microtype}
\DeclareMicrotypeSet*[tracking]{my}%
  { font = */*/*/sc/* }%
\SetTracking{ encoding = *, shape = sc }{ 45 }

\title[Antipodal Sets and Topology]{Maximal Antipodal Sets and the\\ Topology of Generalised Symmetric Spaces}
\def\titl{Maximal Antipodal Sets and the\\ Topology of Generalised Symmetric Spaces}
\def\auth{Manuel Amann}
\date{January 5th, 2021}

\thanks{}

\author{\auth}

\usepackage[english]{babel}

\usepackage[colorlinks,pdftex, plainpages=false]{hyperref}
\hypersetup{pdftitle=\titl, pdfauthor=\auth, pdftoolbar=false,
plainpages=false, hyperindex=true, pdfdisplaydoctitle=true}

\hypersetup{colorlinks,%
citecolor=black,%
filecolor=black,%
linkcolor=black,%
urlcolor=black,%
pdftex}

\usepackage{color}

\usepackage{amsmath, amssymb, amscd, amsthm}
\usepackage{stmaryrd}
\usepackage{latexsym}
\usepackage[all]{xy}
\usepackage{pb-diagram}
\usepackage{rotating}
\usepackage{multicol}
\usepackage{lscape}
\usepackage{wasysym}
\usepackage{longtable}
\usepackage{enumerate}

\xyoption{all}

\newtheorem{theo}{Theorem}[section]
\newtheorem{main}{Theorem}

\newtheorem*{main*}{Theorem}

\newtheorem*{mainprop*}{Proposition}

\newtheorem{mainconj}{Conjecture}

\newtheorem{prop}[theo]{Proposition}
\newtheorem{defi2}[theo]{Definition}
\newtheorem*{defi2*}{Definition}
\newenvironment{defi}{\begin{defi2}\normalfont}{\end{defi2}}
\newenvironment{defi*}{\begin{defi2*}\normalfont}{\end{defi2*}}

\newenvironment{defin*}[1]{\begin{defi2*}[#1]\normalfont}{\end{defi2*}}
\newtheorem*{rem2*}{Remark}
\newenvironment{rem*}{\begin{rem2*}\normalfont}{\hfill$\boxbox$\end{rem2*}}
\newtheorem{rem2}[theo]{Remark}
\newenvironment{rem}{\begin{rem2}\normalfont}{\hfill$\boxbox$\end{rem2}}

\newtheorem{lemma}[theo]{Lemma}
\newtheorem{cor}[theo]{Corollary}
\newtheorem*{cor*}{Corollary}

\newtheorem*{conj*}{Conjecture}
\newtheorem*{theo*}{Theorem}
\newtheorem*{ques*}{Question}
\newtheorem*{mi2}{Main Idea}

\newtheorem{ex2}[theo]{Example}
\newenvironment{ex}{\begin{ex2}\normalfont}{\hfill$\boxbox$\end{ex2}}

\newtheorem{exer2}[theo]{Exercise}

\newtheorem{alg2}[theo]{Algorithm}

\newcommand{\qq}{{\mathbb{Q}}}                                     
\newcommand{\rr}{{\mathbb{R}}}                                     
\newcommand{\s}{{\mathbb{S}}}                                      
\newcommand{\zz}{{\mathbb{Z}}}                                     
\newcommand{\SU}{{\mathbf{SU}}}                                    
\newcommand{\E}{{\mathbf{E}}}                                      
\newcommand{\F}{{\mathbf{F}}}                                      
\newcommand{\G}{{\mathbf{G}}}                                      
\newcommand{\dif} {{\operatorname{d}}}                             
\newcommand{\In} {{\,\subseteq\,}}                                 
\newcommand{\Tor}{{\operatorname{Tor}}}                            
\newcommand{\id}{{\operatorname{id}}}                              
\newcommand{\rk}{{\operatorname{rk\,}}}                            
\newcommand{\Ad}{{\operatorname{Ad}}}                              
\newcommand{\co}{\colon\thinspace}                                 

\newcommand{\comment}[1]{}                                         

\newcommand{\ack}{\noindent\textbf{Acknowledgements. }}            
\newcommand{\str}{\noindent\textbf{Structure of the article. }}    

\newenvironment{prf}{\begin{proof}[\textsc{Proof}]} {\end{proof}}     

\begin{document}

\maketitle \thispagestyle{empty}


\begin{abstract}
We prove several long-standing conjectures by Chen--Nagano on cohomological descriptions of the cardinalities of  maximal antipodal sets in symmetric spaces. We actually extend these conjectures to the setting of generalised symmetric spaces of finite abelian $p$-groups and verify them mostly in this broader context drawing upon techniques from equivariant cohomology theory.
\end{abstract}


\section*{Introduction}

Symmetric spaces have spurred a multitude of research programmes for several decades now. Yet still, ``they and their alikes'' withhold several mysteries to themselves. This may seem surprising in view of the existence of a beautiful classification due to Cartan. However, it becomes more transparent when thinking of their various generalisations or of the lack of knowledge on all sorts of different submanifolds they can harbour.

Not only do we combine both these research directions in this article, we also connect their geometry and topology. More concretely,
\begin{itemize}
\item
the broader class of objects we study (strictly containing symmetric spaces) is given by so-called \emph{generalised symmetric spaces} (more precisely, in the form of \emph{generalised s-structures} for arbitrary finite abelian $p$-groups $\Gamma$ ($p$ a prime number) instead of $\Gamma=\zz_2$ in the symmetric setting, see Definition \ref{def01}),
\item
the submanifolds we consider are so-called \emph{antipodal sets} (as for example a set of antipodal points in a round sphere, see Definition \ref{def02}),
\item
and the connection of geometry and topology will be established by linking cardinalities of (maximal) antipodal sets to cohomological quantities of the respective ambient spaces using equivariant cohomology theory.
\end{itemize}
We are particularly interested in characterising the so-called \emph{antipodal number} $\#_2(M)$ of a symmetric space $M$, respectively its analogue, $\#_\Gamma(M)$, of a generalised symmetric space $M$ (see Section \ref{sec00})---i.e., the supremum over the cardinalities of all (maximal) antipodal sets---via cohomological invariants of $M$ like Euler characteristics and sums of Betti numbers with respective coefficients.

The outcome of these considerations are extensions and verifications of conjectures by Chen--Nagano, which originally were stated for symmetric spaces, mostly in the context of these generalised symmetric spaces. To the best of our knowledge we could trace these conjectures back to the article \cite[p.~53]{Che87} from 1987 (cf.~\cite{CN88}). Let us summarise these three separate conjectures and the one additional question (see \cite[Conjecture 1, p.~27]{Che18}, \cite[Conjecture 2, p.~27]{Che18}, \cite[Conjecture 3, p.~27]{Che18}, \cite[Problem 13.1, p.~28]{Che18}) in a condensed form.
\begin{conj*}[Chen--Nagano, Nagano, 1987] Let $M$ be a connected compact symmetric space. It then holds that
\begin{enumerate}
\item
\begin{align*}
\chi(M) \equiv \#_2(M) \mod 2
\end{align*}
\item and
\begin{align*}
\#_2 M=\dim H_*(M;\zz_2)
\end{align*}
\item
$\#_2 M$ is the smallest number of cells needed for a
CW-structure of $M$.
\item
$\#_2 M> \dim H_*(M;\rr)$ implies that the integral homology of $M$ has $2$-torsion.
\end{enumerate}
\end{conj*}
The equality $\#_2(M)=\chi(M)$ holds on many symmetric spaces and fails for many others (see \cite[Remark 7.3, p.~51]{Che87}).
The second conjecture was verified on symmetric R-spaces by Takeuchi (see \cite{Tak89}). In \cite[p.~54]{Che87} it is stated that one may check $\#_2 M\geq \dim H_*(M;\rr)$ for $M$ simply-connected compact symmetric. There is indeed older knowledge (see \cite{CN88} especially for simply-connected spaces) and recent progress on explicitly classifying maximal antipodal sets of symmetric spaces (for example see \cite{Yu18}, \cite{TT20}). On the one hand it seems these classifications are not yet complete, on the other hand, our arguments in tackling with these conjectures are much shorter and also more general than an explicit classification attempt to maximal antipodal sets in symmetric spaces.

\bigskip

This article proves these conjectures mostly in the much broader context of generalised symmetric spaces with $\Gamma$ a finite abelian $p$-group ($p$ prime)---the only exception being ``$\geq$'' in (2), which we cannot establish, not even for symmetric spaces (see Remark \ref{rem01}). The full scope of Part (3) of the conjecture is not clear to us, yet we show that in general a minimal number of cells may \emph{strictly exceed} the cardinality of any antipodal set (see Example \ref{ex01}) whilst the respective inequality always holds.

We recall that the class of generalised symmetric spaces extends the class of symmetric spaces by many concrete examples. It is also a very interesting task to understand which (topological) properties are shared by symmetric spaces and their generalisations; and the article can be considered as another contribution to this topic. We discuss all this in Remark \ref{rem02}.

The first theorem addresses Aspects (2) and (4) of the conjecture.
\begin{main}\label{theo01}
Let $M$ be a closed connected manifold with a
generalised s-structure for a finite abelian $p$-group $\Gamma$ ($p$ prime).
We have that
\begin{align*}
\#_\Gamma M\leq \dim H^*(M;\zz_p)
\end{align*}
In particular, if $\#_\Gamma M> \dim H_*(M;\rr)$, then $H_*(M;\zz)$ has $p$-torsion.
\end{main}
\begin{rem}\label{rem01}
It seems that from our point of view equality in Theorem \ref{theo01} is related to understanding equivariant formality of the isotropy action with $\zz_p$-coefficients. For rational confirmations of this sort see \cite{Goe12}, \cite{GN16}, \cite{Nos18}, and \cite{AK20}. Combining Theorem \ref{theo01} with Theorem \ref{theo05} below lets us derive Equality in Part (2) of the conjecture (and in the corresponding version for generalised symmetric spaces) at least for several special cases (see Remark \ref{rem03}).
\end{rem}

The next theorem and example deals with the minimal number of cells in a CW-structure, i.e.~with Item (3) in the conjecture.
\begin{main}\label{theo06}
Let $M$ be a closed connected $\Gamma$-symmetric space with $\Gamma$ a finite abelian $p$-group ($p$ prime).
\begin{itemize}
\item
The number of cells needed in a CW-structure of $M$ is bounded from below by $\#_\Gamma(M)$.
\item
If $H^*(M;\zz)$ has torsion other than $p$-torsion, then the number of cells required for a CW-structure of $M$ is strictly larger than $\#_{\Gamma} M$.
\item
Suppose that $\pi_1(M)=0$. If $M$ possesses an antipodal set $A$ with $|A|=\dim H^*(M;\zz_p)$, and if $H^*(M;\zz)$ only has $p$-torsion, then $M$ has a CW-structure with $|A|$ many cells.
\end{itemize}
\end{main}
\begin{ex}\label{ex01}
The cohomology of Riemannian symmetric spaces is well-known (for example see \cite[Section 6.3, p.~147]{MT91}, \cite{Bor53}, \cite{IT77}, etc.) Indeed, most of the spaces do satisfy that their only torsion is $2$-torsion (for example see \cite[Theorem 5.1, p.~241]{IT77} for the exceptional spaces of type I), and in view of the known cardinalities of antipodal sets the theorem does imply the third statement in the conjecture above (in the form that the lower bound is actually even realised) in many cases.

However, there are Riemannian symmetric spaces with torsion other than $2$-torsion. First of all, from \cite[Table 2, p.~524]{Sek88} we cite that there are Riemannian symmetric spaces $G/H$ with
\begin{align*}
(\mathfrak{g}, \mathfrak{h}) \in \{ (\mathfrak{e}_6, \mathfrak{sp}(4)), (\mathfrak{e}_6, \mathfrak{f}(4)), (\mathfrak{su}(p+q), \mathfrak{su}(p)\oplus \mathfrak{su}(q)\oplus \mathfrak{so}(2)) \}
\end{align*}
with fundamental groups $\zz_3$, $\zz_3$, and $\zz_{\gcd(p,q)}$ respectively. This equals their first homology group and produces a respective torsion subgroup in second cohomology via universal coefficients. Hence, due to Theorem \ref{theo06} any of their CW-structures requires strictly more cells than the cardinality of any respective antipodal set (in the last case provided that $\gcd(p,q)$ is not a power of $2$).

Furthermore, even already on the level of simply-connected exceptional Lie groups we find counter-examples. From \cite[Th\'eor\`eme 2.5, p.~223]{Bor61} (cf.~\cite[Theorem 9, p.~1146]{Bor53}, \cite[(3.1), p.~263]{TW74}) we can recall that the cohomology of the groups $\F_4$, $\E_6$, $\E_7$, and $\E_8$ (beside $2$-torsion) has $3$-torsion, and $\E_8$ moreover has $5$-torsion. Hence again, we actually need more cells than the cardinality of any antipodal set.
\end{ex}

Using that the $\Gamma$-structure for a Riemannian symmetric space is explicitly given by geodesic reflections, we can prove Part (1) of the conjecture.
\begin{main}\label{theo02}
Let $M$ be a closed connected Riemannian symmetric space. 
Then it holds that
\begin{align*}
\#_2 M\equiv \chi(M) \mod 2
\end{align*}
More precisely, actually any maximal antipodal set $A$ satisfies $|A|\equiv \chi(M) \mod 2$.
\end{main}
Due to Theorem \ref{theo01} the antipodal number actually is finite.
See Corollaries \ref{cor02} and \ref{cor03} for generalisations to arbitrary $\Gamma$-symmetric spaces.

\bigskip

Let us finally extend a characterisation known for symmetric spaces to our generalised symmetric spaces. In \cite[Theorem 7.1, p.~49]{Che87} it is proved that on a compact symmetric space it holds that
\begin{align*}
\chi(M)\leq \#_2(M)
\end{align*}
We can generalise both the result and (a simplified version of) the proof to generalised symmetric spaces---once having established some more technicalities in this context. Generalised symmetric spaces with $\Gamma$ abelian are known to be homogeneous spaces $M=G/H$ with $G$, $H$ finite dimensional Lie groups (see Section \ref{sec00}).

In the next theorem we assume $G$ and $H$ to be compact, which is guaranteed if the s-structure is Riemannian.
\begin{main}\label{theo05}
Let $M=G/H$ be a closed connected Riemannian $\Gamma$-symmetric space with $\Gamma$ a finite abelian $p$-group ($p$ prime), and with compact (possibly not connected) $G$, $H$. It holds that
\begin{align*}
\chi(M)\leq \#_\Gamma M
\end{align*}
\end{main}
\begin{rem}\label{rem03}
First recall that if $\chi(G/H)>0$, then $\rk G=\rk H$ and
\begin{align*}
\chi(G/H)=\dim H^*(G/H;\qq)
\end{align*}
for compact connected $G$,$H$, since in this case rational cohomology is concentrated in even degrees (see \cite[Theorem 2, p.~217]{Oni94}, or \cite[Example 1, p.~218, and Proposition 32.10, p.~444]{FHT01}). (Note that with $G$ already the component of the identity $G_0\In G$ acts transitively on $M$, and then, for example, $\pi_1(M)=0$ implies that the corresponding stabiliser $H\In G_0$ is connected as well). Let us focus on this case.

Combining Theorems \ref{theo01} and \ref{theo05} it follows that
\begin{align*}
\dim H^*(M;\qq)\leq \#_\Gamma M\leq \dim H^*(M;\zz_p)
\end{align*}
with equalities if and only if $H^*(M;\zz)$ has $p$-torsion (see the second statement of Theorem \ref{theo01}).

Hence in this case we can deduce that
\begin{align*}
\#_\Gamma M=\dim H^*(M;\zz_p)
\end{align*}
(Just to be clear, if there is $p$-torsion this equality might still hold, we just cannot conclude it this way.) In particular, this then confirms Aspect (2) of the conjecture in such a case. For example, this holds if the cohomology of $M$ is torsion-free as amongst others is the case for the symmetric spaces given by the complex Grassmanians---they satisfy the initial prerequisites as well.

Moreover, we can also easily provide generalised symmetric examples satisfying this equality---see Example \ref{ex02}. That is, we can exert explicit computations of antipodal numbers this way in the generalised symmetric setting as well.
\end{rem}

\bigskip

The strategy behind the proofs of these theorems is to draw upon techniques from equivariant cohomology theory applied to a finite abelian $p$-group which we construct using the group $\Gamma$ and an antipodal set---see Section \ref{sec01}. We then relate the antipodal set to the fixed point set of this group acting on $M$. Localisation results from equivariant cohomology connect the $\zz_p$-cohomology of such a fixed point set to the one of $M$. Combining these observations links the cardinality of a (maximal) antipodal set to cohomological invariants of $M$.

We want to advertise rigorously using these equivariant techniques, as they allow for many concise and illuminating arguments for proving the conjectures in much larger generality than, as it seems as of now, otherwise only full classification results possibly could provide.

\bigskip

\str In Section \ref{sec00} we first recall the definitions of generalised s-structures, of generalised symmetric spaces, and of their antipodal sets. Moreover, we recall the necessary concepts from topology, especially from equivariant cohomology theory. In Section \ref{sec01} we provide the construction of the finite abelian $p$-group with respect to which we consider equivariant cohomology, and we prove Theorems \ref{theo01} and \ref{theo06}. Section \ref{sec02} is devoted to the proof of Theorem \ref{theo02}. There, we also provide certain generalisations to $\Gamma$-symmetric spaces. We finish by generalising some structure theory from symmetric spaces to generalised ones, which then permits to prove Theorem \ref{theo05} in Section \ref{sec04}.

\ack The author is grateful to Takashi Sakai for his feedback on previous versions of the article and to Peter Quast for several helpful discussions on the nature of generalised symmetric spaces (and, in particular, for Proposition \ref{prop01}).

The author was supported both by a Heisenberg grant and his research grant AM 342/4-1 of the German Research Foundation; he is moreover a member of the DFG Priority Programme 2026.


\section{Preliminaries}\label{sec00}

Throughout this article, if not stated differently, $\Gamma$ denotes a finite abelian $p$-group (with $p\geq 2$ prime).

\subsection{Generalised symmetric spaces and antipodal sets}

Symmetric spaces are marked by an involutive symmetry around each point $x$. In the case of a Riemannian symmetric space this corresponds to a geodesic reflection $s_x$---clearly having $x$ as an isolated fixed-point. Following and reflecting these geodesics it is then clear that
\begin{align*}
s_x \circ s_y \circ s_x^{-1}=s_{s_x(y)}
\end{align*}
which basically states that the $s_x$ themselves are morphisms compatible with the symmetric structure they establish.
This motivates the following generalisation.
\begin{defi}\label{def01}
A \emph{generalised (regular) s-structure} on a closed manifold $M$ is given by a group $\Gamma$ and a family of group morphisms $(\phi_x)_{x\in M}\co \Gamma\to \operatorname{Diff}(M)$ such that
\begin{enumerate}
\item $x$ is an isolated fixed-point in $M^{\phi_x(\Gamma)}$, and
\item for all $x,y \in M$ and for all $\alpha, \beta\in \Gamma$ it holds that
\begin{align}
\label{eqn01}
\phi_x(\alpha)\circ \phi_y(\beta)\circ \phi_x(\alpha)^{-1}=\phi_{\phi_x(\alpha)(y)} (\alpha\cdot \beta \cdot \alpha^{-1})
\end{align}
\end{enumerate}
We call any such space \emph{Riemannian} if the actions are actually by isometries.

A pseudo-Riemannian manifold which admits such a generalised s-structure is called a \emph{generalised symmetric space} or a \emph{$\Gamma$-symmetric spaces}. If the s-structure is Riemannian, so we call the generalised symmetric space.
\end{defi}
From \cite[Proposition, p.~56]{Lut81} we recall that such a generalised symmetric space $M$ with $\Gamma$ abelian is actually homogeneous $M=G/H$. Due to \cite[Proposition 2.1, p.~56]{Lut81} the automorphism group of the generalised s-structure is a finite dimensional Lie group---hence, without restriction, so are $G$, $H$.
Since the isometry group of a compact manifold is a compact Lie group, in the case of a Riemannian generalised symmetric space $G$ and $H$ can be taken to be compact Lie groups.

Note that for $\Gamma=\zz_2$ acting by isometries this reproduces the definition of a symmetric space; for $\Gamma=\zz_k$ this defines a \emph{regular $k$-symmetric space}.

The object of our study will be generalised symmetric spaces for $\Gamma$ a finite abelian $p$-group, with $p\geq 2$ a prime number.

\begin{rem}\label{rem02}
Generalised symmetric spaces have undergone long research. In particular, this class of spaces extends the class of symmetric spaces by several examples. For example, a classification of $(\zz_2\oplus \zz_2)$-symmetric spaces was obtained in \cite[Theorem 14; Tables 1, 2, 3, 4]{BG08} for the classical cases and \cite[Theorems 1.2, 1.3; Table 1]{Kol09} for exceptional Lie groups~$G$. Respectively, see \cite{KT03} for many examples of $\zz_k$-symmetric spaces.

It is in general a very interesting question to understand which (topological) properties are shared by symmetric spaces and generalised symmetric spaces. For example, they are all \emph{formal} (see \cite{KT03}, \cite{Ste02}, \cite{GN16}) with symmetric spaces also being \emph{geometrically formal} in contrast to the class of $\zz_k$-symmetric spaces in which we can find many examples (like the $6$-symmetric space $\G_2/T^2$) not sharing this property (see \cite{KT03}).

Another such question is the equivariant formality (with rational coefficients) of the isotropy action, which we already commented upon in Remark \ref{rem01}. It is known that this property holds for many such generalised symmetric spaces (see \cite[Question 0.1, p.~2]{AK20}).

One may understand this article also from this perspective of comparing the topological properties of symmetric spaces to their generalised fellows, and we see that most of the properties we investigate are shared by symmetric spaces and generalised ones.
\end{rem}

Recall that with a generalised s-structure on
$M$ one makes the following
\begin{defi}\label{def02}
An \emph{antipodal set} $A\In M$ is a set of elements $x\in M$ such that
\begin{align*}
\phi_x(\alpha)(y)=y \qquad \textrm{for all } x,y\in A \textrm{ and all } \alpha \in \Gamma
\end{align*}
We call an antipodal set $A$ \emph{maximal} if for any antipodal set $A'\In M$ with $A\In A'$ we derive equality $A=A'$.
\end{defi}
Note that from the first property in the definition of an s-structure it follows that an antipodal set is discrete and finite on the compact manifold $M$.

Recall the \emph{antipodal number}
\begin{align*}
\#_{\Gamma} M=\sup\{|A|\mid A \textrm{ is a maximal antipodal set of } (M,g) \}
\end{align*}
which becomes the \emph{$k$-number} in the case of $\Gamma=\zz_k$.

\subsection{Equivariant cohomology theory}

Let us recall some key properties following from localisation results in equivariant cohomology theory. Given a group action of a group $G$ on a space $X$ we denote by $X^G$ its fixed-point set. In this Section we collect several classical theorems connecting the cohomologies of $X$ and of its fixed-point set $X^G$.

Before we do so, let us recall some preliminary easy observations. Cohomology will denote singular cohomology, and by $\dim H^*(X;\Bbbk)$ we denote the \emph{total dimension}, i.e.~the sum over all Betti numbers in the respective coefficient field $\Bbbk\in \{\qq,\rr,\zz_p\}$. Note that for a prime $p\geq 2$ it holds that
\begin{align*}
H^i(X;\zz_p)\cong H_i(X;\zz_p)
\end{align*}
due to the universal coefficient theorem in cohomology/homology for field coefficients (see \cite[Theorem 3.2, p.~195]{Hat02}). Hence we can compute total dimensions of cohomology, which then equal the corresponding total dimensions in homology.

From the very same universal coefficient theorem with $\zz$-coefficients we derive that $H^*(X;\zz)$ has $p$-torsion if and only if so does $H_*(X;\zz)$.

Likewise using universal coefficients, one easily proves that Euler characteristics can equally be computed from homology or cohomology. Moreover, Euler characteristics do not depend on the coefficients.

\bigskip

Let us now finally come to some properties from equivariant cohomology theory. From \cite[(1.5)(b), p.~86]{AP93} or \cite[Theorem 7.10, p.~145]{Bre72} (cf.~\cite[Corollary 1.3.8, p.~30]{AP93} for $p=2$) for positive characteristic as well as from \cite[Corollary 3.1.13, p.~139]{AP93} for characteristic $0$ we recall
\begin{theo}\label{theo03}
Let $X$ be a finite $G$-CW complex with $G$ a compact Lie group.
\begin{itemize}
\item
If $G$ be a finite $p$-group with prime $p\geq 2$, then
\begin{align*}
\chi(X^G)\equiv \chi(X) \mod p
\end{align*}
\item If $G$ is a compact torus $T=\s^1\times \ldots \times \s^1$, then
\begin{align*}
\chi(X^G)=\chi(X)
\end{align*}
\end{itemize}
\end{theo}
As we recall above, the Euler characteristic does not depend on the coefficient field, which explains the similar behaviour in the theorem, although the corresponding localisation results use these different coefficients.

From the slightly more general \cite[Corollary 1.3.8, p.~30]{AP93} for $p=2$, $\Gamma=\zz_2$, and \cite[Corollary 1.4.8, p.~71]{AP93} for $\Gamma=\zz_p$, $p>2$) we cite
\begin{theo}\label{theo04}
Let $X$ be a finite $G$-CW complex, and let $G=\zz_p$ with prime $p\geq 2$. It holds that
\begin{align*}
\dim H^*(X^G;\zz_p)\leq \dim H^*(X;\zz_p)
\end{align*}
\end{theo}
Clearly, this theorem again has a counterpart in characteristic $0$, which we shall not need, however. In order to apply these theorems to smooth actions on smooth manifolds we draw on
\begin{theo}[\protect{\cite[Corollary 7.2, p.~500]{Ill83}}]
Let $G$ be a compact Lie group acting smoothly on a smooth manifold $M$ (with or without boundary). Then $M$ can
be given a $G$-CW-complex structure.
\end{theo}

\begin{rem}\label{rem04}
Theorem \ref{theo03} is formulated for finite $p$-groups (and $\s^1$-tori) already. This results from an inductive reasoning on $\zz_p$-group actions. Also Theorem \ref{theo04} will be applied iteratively. We shall always be interested in the action of a finite abelian $p$-group $T$ for prime $p\geq 2$. In this case due to the Sylow theorems $T$ possesses a subgroup $\zz_p$. Since $T$ is abelian, we have a quotient group $T/\zz_p$ which, obviously, is again a finite abelian $p$-group. Given an effective action of $T$ on a topological space $X$, since $T$ is abelian, it is clear that this induces an action of $T/\zz_p\curvearrowright M^{\zz_p}$ on the $\zz_p$-fixed-point set. We continue in this fashion iteratively considering $\zz_p$-subgroups and iteratively applying Theorem \ref{theo04} in order to deduce that
\begin{align*}
\dim H^*(X;\zz_p)&\geq \dim H^*(X^{\zz_p};\zz_p)
\\&\geq \ldots
\\&\geq  \dim H^*((\ldots(X^{\zz_p})^{\ldots})^{\zz_p};\zz_p)
\\&=\dim H^*(X^{T};\zz_p)
\end{align*}
\end{rem}


\section{Proving Theorems \ref{theo01} and \ref{theo06}}\label{sec01}

Let $M$ be a closed manifold with a generalised s-structure for a finite abelian $p$-group $\Gamma$  (with $p$ prime) of $r$ many summands. Then we fix generators $[1]_i$ of the respective $r$ many cyclic summands of the form $\zz_{p^l}$ for $l\geq 1$. Consider $s_{x,i}:=\phi_x([1]_i)$ with $x$ parametrised over an antipodal set $A\neq \emptyset$---clearly, there is always a non-trivial antipodal set just given by one point. Let $T$ denote the subgroup of the diffeomorphism group generated by these $s_{x,i}$. Let us first make some easy observations using this construction of $T$.
\begin{lemma}\label{lemma01}
Let $A$ be an antipodal set in the generalised symmetric space $M$ with $\Gamma$ a finite abelian $p$-group ($p$ prime). Let $x \in A$ and $y\in M^T$. Then
\begin{align*}
\phi_x(\alpha)\circ \phi_y(\beta)=\phi_y(\beta)\circ \phi_x(\alpha)
\end{align*}
for all $\alpha, \beta\in \Gamma$.
\end{lemma}
\begin{prf} Since $x\in A$ and $y\in M^T$, we have that $\phi_x(\alpha)(y)=y$ for all $\alpha\in \Gamma$, it follows that the defining condition \eqref{eqn01},
\begin{align*}
\phi_x(\alpha)\circ \phi_y(\beta)\circ \phi_x(\alpha)^{-1}=\phi_{\phi_x(\alpha)(y)} (\alpha\cdot \beta \cdot \alpha^{-1})
\end{align*}
together with the commutativity of $\Gamma$ yield
\begin{align*}
\phi_x(\alpha)\circ \phi_y(\beta)\circ \phi_x(\alpha)^{-1}&=\phi_{y} (\alpha\cdot \beta \cdot \alpha^{-1}) \intertext{and}
\phi_x(\alpha)\circ \phi_y(\beta)&=\phi_{y} (\beta) \circ \phi_x(\alpha)
\end{align*}
\end{prf}
Since, by definition of an antipodal set, $A$ is fixed by all $\phi_x(\Gamma)$ for $x\in A$, it is also fixed by the group they generate, i.e.~$\emptyset\neq A\In M^T$. In view of Lemma \ref{lemma01} and of the finiteness of the antipodal set this gives a list of mutually commuting automorphisms $\{s_i\}_{1\leq i\leq n}$ of order a power of $p$ generating $T$. In other words, by the lemma (and, without restriction, not reducing the number of generators) we have a well-defined effective action of a finite abelian $p$-group
\begin{align}\label{eqn02}
T:=T^{\Gamma,A}:=\langle s_1\rangle  \oplus \ldots \oplus \langle s_n\rangle
\end{align}
on $M$. Clearly, $T$ depends on $\Gamma$ and $A$; we suppress this in the notation, however.

\begin{proof}[\textsc{Proof of Theorem \ref{theo01}}]
Denote by $A=\{x_1,x_2,\ldots,x_n\}$ a maximal antipodal set---recall that such a set is discrete and hence finite by compactness---and by $T$ a correspondingly induced finite abelian $p$-group for the action of $\Gamma$ as constructed above in \eqref{eqn02}.

In the case of a Riemannian symmetric space by Lyusternik--Fet and the existence of a closed geodesic, we have that $n\geq 2$.

In general, we may assume that $n\geq 1$ (actually, $n\geq 2$), since otherwise the inequality from the assertion holds trivially.

By definition, the set $A$ is contained in the intersection of the fixed-point sets of the $s_i$, i.e.
\begin{align*}
A\In M^T=M^{s_{1}}\cap \ldots \cap M^{s_{n}}
\end{align*}
Moreover, $A$ is a set of connected components of $M^T$---i.e.~it is open and closed in $M^T$---since by the first property in the definition of a generalised s-structure, any point $x_i\in A$ in $M^T$ (for which we consider $\phi_{x_i}(\Gamma)\In T$) is isolated.

We apply Theorem \ref{theo04} and the inductive approach outlined in Remark \ref{rem04}. Thus we may iteratively pick $\zz_p$-subgroups and successively apply the estimates from the remark. Applying the remark just on each $p$-primary group $\langle s_i\rangle$, i.e.~using that $\dim H^*(M;\zz_p)\geq \dim H^*(M^{\zz_{p^l}};\zz_p)$ for all $l\geq 1$, we may equivalently illustrate this iterative approach in a maybe slightly more intuitive way as follows.

As a first step, we deduce that
\begin{align*}
\dim H^*(M;\zz_p)\geq \dim H^*(M^{s_1};\zz_p)=\sum_i \dim H^*(F^{s_1}_i;\zz_p)
\end{align*}
where the sum runs over the components, $F^{s_1}_i$, of the fixed-point set, $M^{s_1}$. We iterate this argument as in the remark, since due to Lemma \ref{lemma01} the $s_i$ commute with one another, which implies that $s_j$ acts on the fixed-point set of $s_i$ (for any two $1\leq i,j\leq n$). It follows that
\begin{align*}
\dim H^*(M;\zz_p)&\geq \dim H^*(((M^{s_1})^{s_2}\ldots)^{s_n};\zz_p)\\&=\dim H^*(M^{T};\zz_p)\\&=\sum_i \dim H^*(F^{T}_i;\zz_p)
\end{align*}
where the last sum ranges over the fixed-point components, $\bigsqcup_i F_i^T=M^T$, of $M$ under the $T$-action.
As we observed, the discrete set $A$ is a collection of such fixed-point components $F^T_i$, and $\dim H^*(A;\zz_p)=\dim H^0(A;\zz_p)=|A|$. This lets us deduce that
\begin{align*}
\dim H^*(M;\zz_p)&\geq \sum_i \dim H^*(F^{T}_i;\zz_p) \geq |A|
\end{align*}

As $A$ was chosen arbitrarily among antipodal sets, this estimate passes over to the supremum of all cardinalities of maximal antipodal sets whence we deduce the first result
\begin{align*}
\dim H^*(M;\zz_p) \geq \#_\Gamma M
\end{align*}

\bigskip

The second part of the theorem follows easily from the first one. As in the assertion, we assume that $\#_\Gamma M> \dim H_*(M;\rr)$. Hence, since
\begin{align*}
\dim H^*(M;\zz_p)\geq \#_\Gamma M> \dim H^*(M;\rr)=\dim H_*(M;\rr)
\end{align*}
by universal coefficients (see \cite[Theorem 10, p.~246]{Spa81}),
\begin{align*}
H^*(M;\zz_p)\cong H^*(M;\zz)\otimes \zz_p\oplus \Tor(H^{*+1}(M;\zz),\zz_p)
\end{align*}
this strict inequality can only hold if $H^*(M;\zz)$ has $p$-torsion, since $\zz_k\otimes \zz_p\cong\Tor(\zz_k,\zz_p)\cong\zz_{\gcd(k,p)}$ and $\Tor(\zz,\cdot)=0$. As we observed in Section \ref{sec00}, $H^*(M;\zz)$ has $p$-torsion if and only if so does $H_*(M;\zz)$.
\end{proof}

\bigskip

The proof of Theorem \ref{theo06} makes essential use of the first estimate from Theorem \ref{theo01}.
\begin{proof}[\textsc{Proof of Theorem \ref{theo06}}]
From CW-cohomology we recall that the complex $(C^*,\dif)$ underlying $H^*(M;\zz_p)$ is the complex over $\zz_p$ dual to the one freely generated by the cells. It follows that
\begin{align*}
\dim C^*\geq \dim H^*(C^*,\dif)=\dim H^*(M;\zz_p)\geq \#_\Gamma (M)
\end{align*}
for any CW-structure; i.e.~there are at least $\#_\Gamma(M)$ many cells needed in a CW-structure of $M$.

From this chain of inequalities we also see that equality can only hold if the differential $\dif$ on the chain complex modulo $p$ is trivial; let us assume this. The differential $\dif$ is the differential of the complex with integral coefficients reduced modulo $p$. It follows that on the integral chain complex the images of all differentials are divisible by $p$; in other words, the only torsion of $H^*(M;\zz)$ is $p$-torsion.
Formulated slightly differently, if $H^*(M;\zz)$ has $q$-torsion for $q\neq p$ prime, then the differential reduced modulo $p$, i.e.~the differential $\dif$ of $C^*$, is non-trivial, and the number of cells does not equal the total dimension of $\zz_p$-cohomology.

Conversely, if the differential on the integral complex is divisible by $p$, then modulo $p$ every cell constitutes a non-trivial cohomology class, i.e.~$\dim H^*(M;\zz_p)=\dim C^*=\# \textrm{cells}$.

\bigskip

As a compact manifold $M$ possesses a finite CW-structure.
We suppose that $M$ is simply-connected, and that the only torsion in $H^*(M;\zz)$ is $p$-torsion. As we already recalled in Section \ref{sec00}, by the universal coefficient formula for cohomology/homology (see \cite[Theorem 3.2, p.~195]{Hat02}) this is equivalent to $H_*(M;\zz)$ having only $p$-torsion. Using Whitheads's theorem in the form of \cite[Proposition 4C.1, p.~429]{Hat02} we see that this homology now can be realised by a CW-complex which is cellularly homotopy equivalent to the original one and which is constructed  via the following cells. They come either as generators of $\zz$-summands in integral cohomology or in pairs of generators and ``relators'' such that both yield a $p$-torsion summand. As we see from the construction of the cell complex this means that this complex consists of exactly $\dim H_*(M;\zz_p)=\dim H^*(M;\zz_p)$ many cells, and, due to the cellular homotopy equivalence it does provide a CW-structure for $M$. This finishes the proof.
\end{proof}

As we already mentioned in Example \ref{ex01} there are various examples of symmetric spaces satisfying the conditions from Theorem \ref{theo03} as well as symmetric spaces having $3$- and/or $5$- respectively higher torsion. In view of the theorem the following problem hence seems an interesting object of study.
\begin{ques*}
Which (simply-connected) $\Gamma$-symmetric spaces $M$ with $\Gamma$ a finite abelian $p$-group ($p$ prime) satisfy that the only torsion of $H^*(M;\zz)$ is $p$-torsion?
\end{ques*}


\section{Proving Theorem \ref{theo02} and its Corollaries}\label{sec02}

This section deals with the congruence stated in Theorem \ref{theo02}. We can prove it for symmetric spaces and we provide extensions to $\Gamma$-symmetric spaces.

\begin{proof}[\textsc{Proof of Theorem \ref{theo02}}]
The reasoning builds upon the first steps from the proof of Theorem \ref{theo01} together with the following additional observations.

Let us indeed prove that in the case of a Riemannian symmetric space any maximal antipodal set $A$ is not only a collection of fixed-point components of its induced $T$-action on $M$ (as seen in the proof of Theorem \ref{theo01}, $T$ as usual being constructed from $A$ and $\Gamma$ in \eqref{eqn02}), yet, it actually \emph{equals} the $T$-fixed-point set
\begin{align}\label{eqn03}
M^T=A
\end{align}
The key observation for this is that whenever $s_x(y)=y$, then $s_y(x)=x$ for all $x,y \in M^T$, where, as usual, $s_x$ denotes the geodesic reflection at $x$.

Let us quickly justify the latter: in the case of a Riemannian symmetric space $M$ all geodesics are one-parameter groups. This implies that any geodesic loop is actually a closed geodesic. We deduce that whenever $s_{x}(y)=y$ then also $s_y(x)=x$, since (as geodesics only depend on the initial point and the velocity vector) the geodesic reflections at $x$ and at $y$ are both taken with respect to the same closed geodesic.

Let us now prove \eqref{eqn03}. Assume $F$ to be a (non-empty) component of $M^T$ which does not equal a point from $A$. Pick $y\in F$ and consider the geodesic reflection $s_{y}$ at it. Consequently, since $s_x(y)= y$ for all $x\in A$, we also have that $s_{y}(x)=x$ for all $x\in A$ and that $A\sqcup \{y\}$ is an antipodal set strictly containing $A$ by assumption. This contradicts the maximality of $A$. We recall (see Theorem \ref{theo03}) that
\begin{align*}
\chi(M^T) \equiv \chi(M) \mod 2
\end{align*}
since $T$ is a finite (abelian) $2$-group. As $M^T=A$ is a finite (maximal antipodal) set, i.e., in particular, $|A|=|M^T|=\chi(M^T)$, we can conclude that
\begin{align*}
\chi(M)\equiv  \chi(M^T) \equiv |A|\equiv \#_2(M) \mod 2
\end{align*}
by passing to the supremum of all lengths of maximal antipodal sets. (We recall again from Theorem \ref{theo01} that the antipodal number actually is finite, and the supremum is a maximum.)
\end{proof}

In view of Theorem \ref{theo03} also holding for odd prime numbers, the previous proof yields that Theorem \ref{theo02} directly generalises to all $\Gamma$-symmetric spaces with $\Gamma$ a finite abelian $p$-group ($p$ prime), $A$ a maximal antipodal set with associated finite abelian $p$-group $T$ (see \eqref{eqn02}) provided that
\begin{align}\label{eqn04}
\phi_y(\alpha)(x)=x \qquad \textrm{ for all } y\in M^T, x\in A, \alpha \in \Gamma
\end{align}

\begin{cor}\label{cor02}
Let $M$ be a closed connected $\Gamma$-symmetric space with $\Gamma$ a finite abelian $p$-group ($p$ prime) satisfying \eqref{eqn04}.
We obtain that
\begin{align*}
\#_\Gamma M\equiv \chi(M) \mod p
\end{align*}
More precisely, actually any maximal antipodal set $A$ satisfies $|A|\equiv \chi(M) \mod p$.
\end{cor}
A generalised $\Gamma$-symmetric space $M$ (for a finite abelian group $\Gamma$) is homogeneous (see Section \ref{sec00}). In the next assertion we additionally assume the compactness of $G$ and $H$. When speaking about the finite abelian $p$-group associated to an antipodal set we refer to \eqref{eqn02}, as usual.
\begin{cor}\label{cor03}
Suppose that for some choice of compact Lie groups $G$, $H$ with $M=G/H$ it holds $\phi_x(\Gamma)\In G$ for all $x\in M$ and that for a maximal antipodal set $A$ the associated finite abelian $p$-group $T$ is a maximal abelian $p$-subgroup of $G$, then
\begin{align*}
 |A| \equiv \chi(M) \mod p
\end{align*}
\end{cor}
\begin{prf}
We prove that Property \eqref{eqn04} holds in this case whence Corollary \ref{cor02} applies.

Suppose that there exists $y\in M^T\setminus A$. We show \eqref{eqn04}, namely that $\phi_y(\alpha)(x)=x$ for all $x\in A, \alpha \in \Gamma$, which contradicts the maximality of $A$, and thus the existence of $y$ whence $A=M^T$. Hence the arguments of the proof of Theorem \ref{theo02} apply verbatim respectively ``modulo $p$'' instead of ``modulo $2$''.

Since $y\in M^T$, it holds that
\begin{align*}
\phi_y(\alpha)(x)=\phi_{k(y)}(\alpha)(x)=k^{-1}\cdot\phi_y(\alpha)(x)\cdot k
\end{align*}
for all $k\in T$. It follows that $\phi_y(\Gamma)$ commutes with $T$. Thus $T':=\langle T,\phi_y(\Gamma)\rangle\In G$ is a finite abelian $p$-group containing $T$. By the maximality of $T$ in $G$ it follows that $T'=T$, that actually $\phi_y(\Gamma)\In T$, and that \eqref{eqn04} holds.
\end{prf}
It seems interesting to better understand those $\Gamma$-symmetric spaces which satisfy this maximality condition with respect to many maximal antipodal sets.

As for the conidtion that $\phi_x(\Gamma)\In G$ for all $x\in M$ we point the reader to the next Section where refinements of this are discussed at large.

\section{Proving Theorem \ref{theo05}}\label{sec04}

First, we shall need the following generalisation of \cite[Theorem 5.6, p.~424]{Hel78} explained to us by Peter Quast.

Recall that an inner automorphism of a Lie algebra is induced by conjugation. By $\operatorname{fix}(\mathfrak{g},\theta)$ we denote the fixed-point set of the action of an automorphism $\theta$ on a Lie algebra $\mathfrak{g}$.
\begin{prop}\label{prop01}
Let $\mathfrak{g}$ be a compact Lie algebra, and let $\theta$ be an automorphism of $\mathfrak{g}$. Then the following properties are equivalent.
\begin{enumerate}
\item $\theta$ is an inner automorphism.
\item $\rk \operatorname{fix}(\mathfrak{g},\theta)=\rk \mathfrak{g}$
\end{enumerate}
\end{prop}
\begin{prf}
The proof proceeds as follows: First we generalise \cite[Proposition 5.3, p.~423]{Hel78}, which is instrumental in the proof of \cite[Theorem 5.6, p.~424]{Hel78}, to compact Lie algebras. In a second step this allows us to apply the Theorem in larger generality thus yielding the equivalence.

\bigskip

For a compact Lie algebra we have that $\mathfrak{g}=\operatorname{center}(\mathfrak{g}) + [\mathfrak{g},\mathfrak{g}]$ with $\mathfrak{g}':=[\mathfrak{g},\mathfrak{g}]$ being semi-simple. The center lies in every maximal abelian subalgebra $\mathfrak{t}$ of $\mathfrak{g}$.

Now let $\mathfrak{t}$ be such a maximally abelian subalgebra. Then there is some $\mathfrak{t}'\In \mathfrak{g}'$ such that $\mathfrak{t}=\operatorname{center}(\mathfrak{g})+\mathfrak{t}'$. Let $\theta$ denote an automorphism of $\mathfrak{g}$ fixing $\mathfrak{t}$ point-wise. Then $\theta$ leaves $\mathfrak{g}'$ invariant, and its restriction $\theta':=\theta|_{\mathfrak{g}'}$ to the commutator leaves $\mathfrak{t}'$ point-wise fixed. According to \cite[Proposition 5.3, p.~423]{Hel78} there exists some $h'\in \mathfrak{t}'$ such that $e^{\operatorname{ad}(h')}|_{\mathfrak{g}'}=\theta'$. Since $\theta|_{\operatorname{center}{(\mathfrak{g})}}=\id|_{\operatorname{center}(\mathfrak{g})}$, the equality $e^{\operatorname{ad}(h')}=\theta$ holds on all of $\mathfrak{g}$. That is, we proved that a morphism $\theta$ fixing a maximal torus is an inner automorphism, $\theta \in \operatorname{Int}({\mathfrak{g}})$.

\bigskip

Now we observe that the proof of \cite[Theorem 5.6, p.~424]{Hel78} does not make use of $\theta$ being involutive. In view of the preceding arguments it carries over to our situation. For the convenience of the reader we quickly reproduce it here.

Let $G$ be any compact connected Lie group with Lie algebra $\mathfrak{g}$.
Suppose first that $\theta\in \operatorname{Int}(\mathfrak{g})$ is an inner automorphism. Then $\theta = \Ad(g)$ for some $g\in G$. Now $g$
lies in a maximal torus $T\In G$, and $\theta$ leaves the Lie algebra $\mathfrak{t}$ of $T$ point-wise
fixed. Hence $\mathfrak{t} \In \operatorname{fix}(\mathfrak{g},\theta)$, and $\rk \operatorname{fix}(\mathfrak{g},\theta) = \rk \mathfrak{g}$.

Conversely, if $\operatorname{fix}(\mathfrak{g},\theta)$ and $\mathfrak{g}$ have the same rank, there exists a subalgebra $\mathfrak{t}\In \operatorname{fix}(\mathfrak{g},\theta)$ which is maximally
abelian in $\mathfrak{g}$. Our considerations above then show that $\theta\in \operatorname{Int}(\mathfrak{g})$.
\end{prf}
By $G_0$ we denote the connected component of the identity of $G$. We again draw on the fact that a connected $\Gamma$-space $M$ (with $\Gamma$ abelian) is homogeneous $M=G/H$ with the transvection group (generated by all the $\phi_x(\alpha)\phi_y(\alpha)^{-1}$ for $x,y\in M$, $\alpha\in \Gamma$) acting transitively on $M$ (see Section \ref{sec00}). We shall always assume $G$ and $H$ to be compact.
From \cite[p.~57]{Lut81} we recall further that the stabiliser group $H$ in $x_0$ is a collection of connected components of the set of elements fixed by all the
\begin{align*}
g\mapsto \tilde \phi(g)(\alpha):=\phi_{x_0}(\alpha)\circ \phi_{g(x_0)}(\alpha)^{-1}\circ g=\phi_{x_0}(\alpha)\circ g\circ \phi_{x_0}(\alpha)^{-1}
\end{align*}
and that the original $\Gamma$-structure transcribes to a $\Gamma$-structure on $G/H$ given as
\begin{align*}
\phi_{|g|}(\alpha)(|g'|)=|g \tilde\phi(\alpha)(g^{-1}g')|
\end{align*}
We refer to this $\Gamma$-structure as the ``homogeneous" one. Note that
\begin{align*}
|g\tilde \phi(\alpha)(g^{-1}g')|&=|g\phi_{x_0}(\alpha) g^{-1}g'\phi_{x_0}(\alpha)^{-1}|
\\&=|(g\phi_{x_0}(\alpha) g^{-1})(g'\phi_{x_0}(\alpha)^{-1}g'^{-1})g'|
\\&=(\phi_x(\alpha)\circ \phi_y^{-1}(\alpha)) (y)
\\&=\phi_x(\alpha)(y)
\end{align*}
whenever $x=gx_0=|g|$, $y=g'x_0=|g'|$ under the identification $M=G/H$.
We can prove some simple properties using this correspondence. For the sake of simplicity we shall suppress the dependence on $\alpha$ in the following
\begin{lemma}\label{lemma02}
Assume that $\rk G=\rk H$ with a common maximal torus $T$. Suppose that $x_0\in M^T$.
\begin{enumerate}
\item
For $x,y\in M$ it holds that
\begin{align*}
\phi_x(y)=y \iff \phi_{|g|}(|g'|)=|g'|
\end{align*}
whenever $x=gx_0=|g|$, $y=g'x_0=|g'|$ under the identification $M=G/H$. That is, in particular, there is a bijection of antipodal sets of the same cardinality.
\item
The lift to $G$ of $\phi_{|e|}(\cdot)$ in the identity $e\in G$ fixes the maximal torus $T$, i.e.~we have automorphisms (parametrised by $\Gamma$)
\begin{align*}
\phi_e\co G\to G \quad  \phi_e(g')= \tilde\phi(g') \quad \textrm{satisfying} \quad k\mapsto \phi_{e}(k)=k
\end{align*}
for $k\in T$. Hence the analogue holds for the respective derivatives $\theta$ with respect to the maximal abelian subalgebra $\mathfrak{t}$ of $\mathfrak{g}$.
\item As a consequence we have that the homogeneous s-structure consists of inner automorphisms $\phi_{|g|}(\cdot)\in G_0$.
\item
The homogeneous $\Gamma$-structure commutes with the maximal torus, i.e.~
\begin{align*}
k\cdot \phi_{|g|}(|g'|)=\phi_{|kg|}(|kg'|)
\end{align*}
for all $g, g'\in G$, $k\in T$.
\end{enumerate}
\end{lemma}
\begin{prf}
As for the first part, suppose first that $\phi_x(y)=y$.  We compute
\begin{align*}
\phi_{|g|}(|g'|)|&=|g \tilde\phi(g^{-1}g')|
\\&=|g \phi_{x_0}\circ \phi_{g^{-1}g'(x_0)}^{-1}\circ g^{-1}g'|
\\&=|(g \circ \phi_{x_0}\circ  g^{-1})g'|
\\&=|\phi_{x}g'|
\\&=|g'|
\end{align*}
Conversely, suppose that $ \phi_{|g|}(|g'|)=|g'|$. It follows that
\begin{align*}
|g'|&=\phi_{|g|}(|g'|)
\\&=|g \tilde\phi (g^{-1}g')|
\\&=|g \circ \phi_{x_0}\circ g^{-1}g'\circ \phi_{x_0}^{-1}|
\\&=|(g \circ \phi_{x_0}\circ g^{-1})\circ(g'\circ \phi_{x_0}^{-1}\circ g'^{-1})\circ g'|
\\&=|\phi_x\circ \phi_y \circ g'|
\\&=|\phi_x g'|
\end{align*}
yielding that $\phi_x y=y$.

\bigskip

We prove the second assertion as follows. First note that, since $x_0\in M^T$ and $k\in T$, we have that
\begin{align}\label{eqn05}
\phi_{x_0}=\phi_{kx_0}=k\cdot\phi_{x_0}\cdot k^{-1}
\end{align}
whence $k$ commutes with $\phi_{x_0}$. It follows that
\begin{align*}
\phi_{e}(k)=\tilde \phi(k)
=\phi_{x_0}\circ k\circ \phi_{x_0}^{-1}
=k
\end{align*}

\bigskip

For the third assertion we first compute that
\begin{align*}
\phi_{|g|}=g\tilde\phi(g^{-1}g')=g\phi_{x_0}\circ g^{-1}g'\circ \phi_{x_0}^{-1}=g\phi_e(g^{-1}g')
\end{align*}
Since $M$ is connected, the connected component of the identity $G_0$ acts transitively upon it. Hence, without restriction, $g\in G_0$, and the diffeomorphism of $G$, $g\phi_e(g^{-1}g')$, is homotopic to $\phi_e(g')$. From Item (2), however, we know that the derivative of the latter fixes a maximal abelian subalgebra $\mathfrak{t}$ of $\mathfrak{g}$. By Proposition \ref{prop01} we hence deduce that $\phi_{e}(\cdot)$ is an inner automorphism of $G$. Hence its projection $\phi_{|e|}(\cdot)$ is in $G_0$, and it is homotopic to $\phi_{|g|}(\cdot)$. Consequently, the $\phi_{|g|}(\cdot)$ for all $|g|\in G/H$ lie in $G_0$.

\bigskip

As four Assertion (4) we compute
\begin{align*}
k\cdot \phi_{|g|}(|g'|)=kg\tilde \phi(g^{-1}g')=kg\tilde \phi((kg)^{-1} kg')=\phi_{kg}(kg')
\end{align*}
\end{prf}

Let us now prove Theorem \ref{theo05}. We start by following along the lines of the proof of \cite[Theorem 7.1, p.~49]{Che87}. However, we shall alter arguments in the course of the proof.
\begin{proof}[\textsc{Proof of Theorem \ref{theo05}}]
As we recalled in Section \ref{sec00}, the manifold $M$ is a homogeneous space $M=G/H$; we assume $G$ and $H$ to be compact. We may assume that $\rk G=\rk H$, since otherwise $\chi(M)=0$, and the assertion is trivially true. We may pick a common maximal torus $T^{\max}\In H\In G$ of $H$ and $G$, which then acts with fixed-points by left multiplication (i.e.~via the isotropy action) on $G/H=M$. It is our (first) goal to identify $M^{T^{\max}}$ as an antipodal set. For this, in view of Lemma \ref{lemma02}.1 we may work with the induced homogeneous s-structure as depicted above.

Let now $y\in M^{T^{\max}}$ be a point from this fixed-point set. In Lemma \ref{lemma02}.4 we have seen that this homogeneous s-structure commutes with the action of $T^{\max}$. In Lemma \ref{lemma02}.3 we observed that the $\phi_{|g|}(\alpha)(\cdot)$ lie in the component of the identity $G_0\In G$.

Combining these two properties from the last paragraph and recalling that a maximal torus is a maximal abelian subgroup (see \cite[Corollary 2.7, p.~287]{Hel78}), we derive that $\phi_y(\Gamma)\In T^{\max}$ for all $y\in M^{T^{\max}}$. It follows that $\phi_y(\Gamma)(x)=x$ for all $x,y\in M^{T^{\max}}$, and $M^{T^{\max}}$ is an antipodal set.

\bigskip
In the second part of this proof we now apply equivariant cohomology theory using this concrete antipodal set.
Theorem \ref{theo03} yields that $\chi(M)=\chi(M^{T^{\max}})$. Passing to the supremum over the lengths of all (maximal) antipodal sets lets us derive the result
\begin{align*}
\chi(M)=\chi(M^{T^{\max}})=|M^{T^{\max}}|\leq \#_\Gamma (M)
\end{align*}
\end{proof}

In Remark \ref{rem03} we already discussed applications of Theorem \ref{theo05} in view of equality in Aspect (2) of the main conjecture, i.e.~in view of actually computing cardinalities of antipodal sets via $p$-cohomology. We promised a generalised symmetric example, namely
\begin{ex}\label{ex02}
From \cite[p.~500]{KT03} we recall that the flag manifold, $\SU(n)/T^{n-1}$ ($n\geq 2$) is a $\zz_{n}$-symmetric space. Choose $n=p$ prime. Clearly, $\rk \SU(n)=\rk T^{n-1}=n-1$, and both groups are connected. We see that the computation of the real cohomology of this space in \cite[Proposition 1, p.~497]{KT03} (for example via a spectral sequence argument) carries over to give the analogue integral cohomology structure, which is torsion-free. Hence our results apply and yield that
\begin{align*}
\#_p (\SU(p)/T^{p-1})&=\dim H^*(\SU(p)/T^{p-1};\zz_p)
\\&=\dim H^*(\SU(p)/T^{p-1};\qq)
\\&=\chi(\SU(p)/T^{p-1})
\\&=|W(\SU(p))|
\\&=2^p
\end{align*}
where $W(\SU(p))=S_p$ denotes the Weyl group of $\SU(p)$, the symmetric group $S_p$. We used that $\rk \SU(p)=\rk T^{p-1}$ together with \cite[Proposition 3.31, p.~127]{FOT08}.
\end{ex}



\def\cprime{$'$}


\vfill

\begin{center}
\noindent
\begin{minipage}{\linewidth}
\small \noindent \textsc
{Manuel Amann} \\
\textsc{Institut f\"ur Mathematik}\\
\textsc{Differentialgeometrie}\\
\textsc{Universit\"at Augsburg}\\
\textsc{Universit\"atsstra\ss{}e 14 }\\
\textsc{86159 Augsburg}\\
\textsc{Germany}\\
[1ex]
\footnotesize
\textsf{manuel.amann@math.uni-augsburg.de}\\
\textsf{www.uni-augsburg.de/de/fakultaet/mntf/math/prof/diff/team/dr-habil-manuel-amann/}
\end{minipage}
\end{center}

\end{document}